\documentclass[twoside,12pt]{amsart}
\usepackage{mathrsfs}
\usepackage{amsfonts}
\newtheorem{thm}{Theorem}[section]

\newtheorem{lem}[thm]{Lemma}
\newtheorem{prop}[thm]{Proposition}

\makeatletter
\def\serieslogo@{}
\makeatother
\makeatletter
\def\@setcopyright{}
\makeatother

\theoremstyle{definition}

\theoremstyle{remark}
\newtheorem{rem}{Remark}[section]

\begin{document}

\title{stability of Steady solutions to Reaction-Hyperbolic Systems for Axonal Transport}

\author{Hao YAN}
\address{Zhou Pei-Yuan Center for Appl. Math.\\
Tsinghua University\\
Beijing 100084, China} \email{yanhao06@mails.tsinghua.edu.cn}

\author{Wen-An YONG}
\address{Zhou Pei-Yuan Center for Appl. Math.\\
Tsinghua University\\
Beijing 100084, China} \email{wayong@tsinghua.edu.cn}




\begin{abstract}
This paper is concerned with the stability of steady solutions to
initial-boundary-value problems of reaction-hyperbolic systems for
axonal transport. Under proper structural assumptions, we clarify
the relaxation structure of the reaction-hyperbolic systems and show
the time-asymptotic stability of steady solutions or relaxation
boundary-layers.

\vskip 10 pt

\end{abstract}

\maketitle \markboth{Hao YAN and W.-A. YONG}{Stability of steady
solutions for Reaction-Hyperbolic Systems for Axonal Transport}

\section{Introduction}
This work is concerned with the following reaction-hyperbolic system
\begin{equation}\label{pdel}
\partial_t u_i + \lambda_i\partial_x u_i = \sum_{j=1}^rk_{ij}u_j, \qquad i =
1, 2, \cdots, r
\end{equation}
on the quarter-plane $x,t\ge0$. Here $u_i=u_i(x, t)(i=1,2,\cdots,r)$
are unknown functions, $\lambda_i(i=1,2,\cdots,r)$ and
$k_{ij}(i,j=1,2,\cdots,r)$ are given constants. It is well believed
\cite{RB} that such models describe the axonal transport in
neuroscience.

The axonal transport is important for the maintenance and functions
of nerve cells or {\it neurons}.
A neuron consists of three parts mainly: cell body, dendrites and a
single axon. The axon is a long and thin pipe whose length can
exceed 10,000 times its diameter. It is this axon that distinguishes
neurons from cells in other organs or tissues. The axon is
responsible for signal transmission in the nervous system. Its
cytoplasm does not contain rough endoplasmic reticulum and therefore
its proteins can only be transported from the cell body, where all
proteins are synthesized. In addition, the axonal transport is also
used by the neuroscientists to trace the connections in the brain.

The transport proceeds as follows. Proteins are stored in vesicles
as cargos. The vesicles are attached to kinesin (anterograde motors)
or dynein (retrograde motors) proteins. These motor proteins drive
the vesicles to walk along the cytoskeletal microtubules as track.
Here the kinesin proteins move the vesicles from the cell body to
synapse (anterograde transport), while the dynein proteins move the
vesicles in the opposite direction (retrograde transport). During
the transport, many biochemical processes are possible. For example,
the cargos can leave its track, can switch its motor proteins from
kinesin to dynein or vice verse, and can move back onto the track.
Thus, we can divide the cargos into a number of subpopulations, such
as free vesicles,
vesicle-kinesin compounds off track, moving vesicle-dynein compounds
on track, etc.

As the axon is long and thin, it is reasonable to assume the
transport only along the longitudinal direction of the axon. Denote
by $x>0$ the distance down the axon from the cell body
which is located at $x = 0$. Let $u_i = u_i(x, t)$ be the
concentration at space-time $(x, t)$ of the $i$-th subpopulations.
According to Reed and Blum \cite{RB}, the mathematical model for
axonal transport is partial differential equations of the form
\begin{equation}\label{pde0}
\partial_t u_i + \lambda_i\partial_x u_i = F_i (u_1, u_2, \cdots, u_r), \qquad i =
1, 2, \cdots, r
\end{equation}
defined on the quarter-plane $x, t\ge0$. Here the term
$\lambda_i\partial_x u_i$ accounts for the transport of the $i$-th
subpopulation with constant velocity $\lambda_i$, and $F_i (u_1,
u_2, \cdots, u_r)$ describes the biochemical processes of the
constituents. It is well recognized that the biochemical processes
are much faster than the transport in biosystems. Thus, it is more
proper to rewrite the general equation $(\ref{pde0})$ as
\begin{equation}\label{pde0*}
\partial_t u_i + \lambda_i\partial_x u_i =\frac{1}{\epsilon} F_i (u_1, u_2, \cdots, u_r),
\qquad i =1, 2, \cdots, r .
\end{equation}
Here the small parameter $\epsilon>0$ characterizes the fact that
the biochemical processes are much faster than the transport.

On the other hand, in a typical experiment for investigating axonal
transport, a large quantity of radiolabeled amino acids are injected
into a nerve ganglion. Thus, the proteins are synthesized with the
injected amino acids in the cell body and transported along the
axon, and finally the radioactivity appears in the axon in a few
hours. The wave of radiolabeled proteins travels at approximately
constant velocity. The shapes and speeds of the profiles depend on
nerves, animals, and proteins being studied. Usually the amino acids
are continuously available in the experiment, and the sharp
approximately parallel the wave fronts which suggest traveling
waves. There are at least two such systems which carry
membrane-bounded organelles and cytoplasm separately. For more
details about the experiments for axonal transport, the reader can
refer to \cite{Can,GB,Oc}.

In order to explain the approximate traveling waves observed in
experiments, the authors of \cite{CBF, FC, FH, RVB} studied the
linear case (\ref{pdel}), where $F_i(u_1, u_2, \cdots, u_r)$ is
linear with respect to the $u_j$'s. A concrete example can be found
in \cite{CBF} on neurofilament transport. Because parabolic
equations typically admit traveling wave solutions, these authors
all related the system (\ref{pde0*}) to a parabolic equation--- the
diffusive limit of (\ref{pde0*}). Especially, in \cite{FH} Friedman
and Hu used parabolic-type estimates to analyse the diffusive limit
of the linear systems.
However, it is known \cite{YZ, DY} that, unlike the conservation
laws, hyperbolic systems with relaxation also allow non-trivial
traveling wave solutions.

In this project, we intend to explain the experimental observation
by directly studying steady solutions to the initial-boundary-value
problems of the above hyperbolic-reaction systems. Steady solutions
very much look like the traveling waves but they should called
(relaxation) boundary-layers due to the presence of the boundary
$x=0$. Because they can be observed experimentally, we believe that
the steady solutions are stable. Therefore, the goal of this paper
is to investigate the time-asymptotic stability of the steady
solutions under proper structural assumptions on the system.

The standard structural assumptions on the axonal transport models
$(\ref{pdel})$ read as
\begin{itemize}
\item[(H1).] $k_{ij}\ge0$ if $i\ne j$;\\

\item[(H2).] (conservation of mass) $\sum_{i=1}^rk_{ij}=0$ for $j=1,2,\cdots ,r$;\\

\item[(H3).] (irreducibility) for any $i_0\ne i_1$, there is a sequence of indices $j_1, j_2,\cdots,j_l$

such that $i_0=j_1, i_1=j_l$ and $k_{j_mj_{m+1}}>0$ for $1\le m\le l-1$;\\

\item[(H4).] there exist $i$ and $j$ such that $\lambda_i\ne \lambda_j$.
\end{itemize}
These assumptions are taken from \cite{FH}. See also
\cite{Ca,CBF,FC,RB,RVB}. Remark that (H4) ensures the
hyperbolic-reaction system $(\ref{pdel})$ won't degenerate into a
system of ordinary differential equations.

As a first step of this project, we consider the linear system
$(\ref{pdel})$ and additionally make the following technical
assumption in this paper that
\begin{itemize}
\item[(H5).] $\lambda_i>0,\qquad i=1,2,\cdots,r$.
\end{itemize}
In the future, we will remove or relax this (H5).

With the above assumptions, we formulate the problem as follows.
Write
$$
U = (u_{1}, u_{2}, \cdots, u_{r})^{T}, \quad \Lambda =
\mbox{diag}(\lambda_1,\lambda_2,\cdots,\lambda_ r), \quad K =
[k_{ij}]_{n\times n}.
$$
Then $(\ref{pdel})$ can be written as
\begin{eqnarray}\label{pde3}
U_t + \Lambda U_x = KU.
\end{eqnarray}
The initial and boundary data are
\begin{eqnarray}\label{value}
U|_{t=0} = U_0(x), \qquad U|_{x=0} = U_0(0) .
\end{eqnarray}
Here the second equality implies the consistency condition which is
necessary for the above initial-boundary-value problems to have a
continuous solution. Moreover, we assume that
\begin{eqnarray}\label{value1}
\Lambda U_{0x}(0) = KU_0(0),
\end{eqnarray}
which is
necessary for the above problems to have a
continuously differentiable solution.
Our main results are the time-asymptotic
stability of steady solutions to $(\ref{pde3})$ together with
$(\ref{value})$.

For references on nonlinear systems for axonal transport, we mention
that Carr in \cite{Ca} discussed the existence of global classical
solutions to a class of nonlinear models. Recently, in \cite{YY} we
clarified the relaxation structure of nonlinear models in \cite{Ca}
and verified the relaxation limit of BV-solutions to the Cauchy
problems.

The rest of this paper is organized as follows. In Section 2 we
discuss the relaxation structure of the hyperbolic-reaction system
(\ref{pdel}). Section 3 is devoted to the existence and stability of
the steady solutions for $(\ref{pde3})$.

\section{Relaxation Structure}
\setcounter{equation}{0}

In this section, we show that the hyperbolic-reaction system
(\ref{pdel}) possesses the relaxation structure formulated in
\cite{Y2}, provided that the assumptions (H1)--(H3) hold. We start
with the following elementary fact as (a) of Lemma 1 proved in
\cite{RVB}.
\begin{lem}\label{lemma1}
Under the assumptions (H1)--(H3), the kernel of the matrix $K$ is
one-dimensional and is spanned by a vector with strictly positive
entries.
\end{lem}
%
On the basis of this fact, we can show
\begin{lem}\label{lemma2}
$0$ is a single eigenvalue of $K$.
\end{lem}
\begin{proof}
Set
\begin{eqnarray*}
L_1=\left( \begin{array}{cccccc}
1 & e_1\\
0 & I_{r-1}
\end{array} \right),
\end{eqnarray*}
where $e_1=(1,1,\cdots,1)$ is a vector in $\mathbb{R}^{r-1}$,
$I_{r-1}$ is the unit matrix of dimension $(r-1)$, and the
superscript T denotes the transpose of vectors or matrices. Then the
inverse of $L_1$ is
\begin{eqnarray*}
L_1^{-1}=\left( \begin{array}{cccccc}
1 & -e_1\\
0 & I_{r-1}
\end{array} \right).
\end{eqnarray*}
Using the partition of $L_1$, we rewrite $K$ as
\begin{eqnarray*}
K=\left( \begin{array}{cccccc}
k_{11} & \alpha\\
\beta & K_1
\end{array} \right),
\end{eqnarray*}
where $\alpha=(k_{12},k_{13},\cdots,k_{1r}),
\beta=(k_{21},k_{31},\cdots,k_{r1})^T$, and $K_1$ is the
$(r-1)\times (r-1)$-matrix. By a direct calculation using the
assumption (H2), we obtain
\begin{eqnarray*}
L_1KL_1^{-1}=\left( \begin{array}{cccccc}
1 & e_1\\
0 & I_{r-1}
\end{array} \right)\left( \begin{array}{cccccc}
k_{11} & \alpha\\
\beta & K_1
\end{array} \right)\left( \begin{array}{cccccc}
1 & -e_1\\
0 & I_{r-1}
\end{array} \right)=\left( \begin{array}{cccccc}
0 & 0\\
\beta & K_2
\end{array} \right)
\end{eqnarray*}
with $K_2 = K_1 - \beta e_1$.

Thus, it reduces to show that the sub-matrix $K_2$ is invertible.
Otherwise, there is a vector $\eta=(\eta_1, \eta_2, \cdots,
\eta_{r-1})^T$ such that $K_2\eta=0.$ Take
\begin{eqnarray*}
\xi=L_1^{-1}\left( \begin{array}{cccccc}
0\\
\eta
\end{array}\right).
\end{eqnarray*}
Then it holds that
\begin{eqnarray*}
K\xi=KL_1^{-1}\left( \begin{array}{cccccc}
0\\
\eta
\end{array}\right)=L_1^{-1}\left( \begin{array}{cccccc}
0 & 0\\
\beta & K_2
\end{array} \right)\left( \begin{array}{cccccc}
0\\
\eta
\end{array}\right)=L_1^{-1}\left( \begin{array}{cccccc}
0\\
K_2\eta
\end{array}\right)=0.
\end{eqnarray*}
This shows that $\xi$ is a vector in the kernel of $K$. On the other
hand, we have
\begin{eqnarray*}
\xi=L_1^{-1}\left( \begin{array}{cccccc}
0\\
\eta
\end{array}\right)=\left( \begin{array}{cccccc}
1 & -e_1\\
0 & I_{r-1}
\end{array} \right)\left( \begin{array}{cccccc}
0\\
\eta_1\\
\eta_2\\
\cdots\\
\eta_{r-1}
\end{array}\right)=\left( \begin{array}{cccccc}
-\eta_1-\eta_2-\cdots-\eta_{r-1}\\
\eta_1\\
\eta_2\\
\cdots\\
\eta_{r-1}
\end{array}\right).
\end{eqnarray*}
Obviously, such a $\xi$ can not be in the kernel of $K$ spanned by a
vector with positive entries. This contradicts Lemma $\ref{lemma1}$.
Therefore, $K_2$ is invertible and $0$ is a single eigenvalue of
$K$. This completes the proof.
\end{proof}
\begin{rem}\label{rem1}
Following the above proof, we have
\begin{eqnarray*}
\left( \begin{array}{cccccc}
1 &0  \\
K_2^{-1}\beta &I_{r-1}
\end{array} \right)L_1K = \left( \begin{array}{cccccc}
0 &0  \\
0 &K_2
\end{array} \right)\left( \begin{array}{cccccc}
1 &0  \\
K_2^{-1}\beta &I_{r-1}
\end{array} \right)L_1.
\end{eqnarray*}
On the other hand, it is not difficult to deduce from the Gershgorin circle theorem that
non-zero eigenvalues of $K$ have negative real parts. Therefore, $K_2$ is stable and (i)
of the first stability condition in \cite{Y2} is verified.
\end{rem}

Furthermore, we have
\begin{lem}\label{lemma3}
Under the assumptions (H1)--(H3), there exist a positive definite
diagonal matrix $A_0$, an orthogonal matrix $P$, and a positive
definite diagonal matrix $S$ such that
\begin{eqnarray*}
&A_0\Lambda=\Lambda A_0, \\
&A_0K+K^TA_0=-P^T\left( \begin{array}{cccccc}
0 &0  \\
0 &S
\end{array} \right)P.
\end{eqnarray*}
Furthermore, the first column of $P^T$ is in the kernel of $K$.
\end{lem}
\begin{proof}
Let $\xi=(\xi_1, \xi_2, \cdots, \xi_r)^T$ is an eigenvector of the
matrix $K$, associated with the eigenvalue $0$. By Lemma
$\ref{lemma1}$, we may assume that $\xi_i>0$ for each $i$. Define
\begin{eqnarray*}
D=\left( \begin{array}{cccccc}
\xi_1 & & & &  \\
& &\xi_2&  &\\
& & &\cdots &\\
& & & &\xi_r
\end{array} \right).
\end{eqnarray*}
It is obvious that the matrix $KD$ satisfies the assumptions
(H1)--(H3) as well and
\begin{eqnarray*}
KD\left( \begin{array}{cccccc}
1\\
1\\
\cdots\\
1
\end{array}\right)=K\left( \begin{array}{cccccc}
\xi_1\\
\xi_2\\
\cdots\\
\xi_r
\end{array}\right)=0.
\end{eqnarray*}
Namely, the sum of each row of $KD$ is also 0. Moreover, it is easy
to see that the symmetric matrix $KD+DK^T$ fulfils the (H1)--(H3),
too.

Take $A_0=D^{-1}$ and it is clear that $A_0$ is positive definite
diagonal and $A_0\Lambda=\Lambda A_0$. Moreover, the symmetric
matrix
\begin{eqnarray*}
A_0K+K^TA_0=D^{-1}(KD+DK^T)D^{-1}
\end{eqnarray*}
also fulfils the assumptions (H1)--H(3). According to Lemmas
$\ref{lemma1}$ and $\ref{lemma2}$, $0$ is a single eigenvalue of
$A_0K+K^TA_0$. Moreover, from the Gershgorin circle theorem it is
not difficult to deduce that non-zero eigenvalues of $A_0K+K^TA_0$
are negative. Thus, there exist an orthogonal matrix $P$ and a
positive definite diagonal matrix $S$ such that
\begin{eqnarray*}
&A_0K+K^TA_0=-P^T\left( \begin{array}{cccccc}
0 &0  \\
0 &S
\end{array} \right)P.
\end{eqnarray*}

Furthermore, since the $r$-vector $(\xi_1,\xi_2,\cdots,\xi_r)^T$ is
in the kernel of $A_0K+K^TA_0$, it is easy to see that the last
$(r-1)$ components of the column vector $P(\xi_1, \xi_2, \cdots,
\xi_r)^T$ are zeros. Namely, the vector $(\xi_1,\xi_2,\cdots,\xi_r)$
is orthogonal to the last $(r-1)$ rows of the orthogonal matrix $P$
and thereby parallels to the first row of $P$. Hence, the first column
of $P^T$ is in the kernel of $K$ and the proof is complete.
\end{proof}

\begin{rem}
Lemma \ref{lemma3}, together with Remark \ref{rem1}, shows that the
reaction-hyperbolic systems satisfying assumptions (H1)--(H3) fulfil the
first stability condition in \cite{Y1, Y2}. However, they do not satisfy
the second stability condition in \cite{Y1, Y2} in general, unless further assumptions
are posed. An important case is that $KD$ is symmetric, which implies the
second stability condition due to Theorem 5.3 in \cite{Y2}.
It is clear that $KD$ is symmetric, provided that the principle of detailed balance
holds (see, e.g., \cite{YGCGM}). For the neurofilament model in \cite{CBF}, $K$ is tri-diagonal
and one can easily see that $KD$ is symmetric. However,
the assumptions (H1)--(H3) do not imply the symmetry of $KD$.
In fact, the $4\times 4$-matrix
\begin{eqnarray*}
K=\left( \begin{array}{cccccc}
-4&1&1&0 \\
2&-3&0&1\\
2&1&-2&0\\
0&1&1&-1
\end{array} \right),
\end{eqnarray*}
satisfies (H1)--(H3) and $(\frac{1}{2},1,1,2)^T$ is an eigenvector associated with 0.
Set $D=diag(\frac{1}{2},1,1,2)$. By a direct calculation, we have
\begin{eqnarray*}
KD=\left( \begin{array}{cccccc}
-2&1&1&0 \\
1&-3&0&2\\
1&1&-2&0\\
0&1&1&-2
\end{array} \right),
\end{eqnarray*}
which is not symmetric.
\end{rem}

But we have
\begin{prop}
For $r\le3$, the matrix $KD$ is symmetric.
\end{prop}
\begin{proof}
For $r=2$, let $K$ be
\begin{eqnarray*}
K=\left( \begin{array}{cccccc}
a &-b  \\
-a &b
\end{array} \right)
\end{eqnarray*}
with $a,b>0$. We take $\xi=(b,a)^T$ and $D=diag(b,a)$. Then
\begin{eqnarray*}
KD=\left( \begin{array}{cccccc}
ab &-ab  \\
-ab &ab
\end{array} \right)
\end{eqnarray*}
is symmetric.

For $r=3$, set
\begin{eqnarray*}
K=\left( \begin{array}{cccccc}
k_{11}& k_{12}&k_{13} \\
k_{21}&k_{22}&k_{23}\\
k_{31}&k_{32}&k_{33}
\end{array} \right)
\end{eqnarray*}
and $D=\mbox{diag}(\xi_1, \xi_2, \xi_3)$.
Moreover, set
$$
a = \xi_2k_{12} - \xi_1k_{21}, b = \xi_3k_{13} - \xi_1k_{31}, c=\xi_3k_{23} - \xi_2k_{32}.
$$
Since
\begin{eqnarray*}
K\left(\begin{array}{cccccc}
\xi_1\\
\xi_2\\
\xi_3
\end{array}\right)=\left( \begin{array}{cccccc}
\xi_1k_{11}+\xi_2k_{12}+\xi_3k_{13} \\
\xi_1k_{21}+\xi_2k_{22}+\xi_3k_{23}\\
\xi_1k_{31}+\xi_2k_{32}+\xi_3k_{33}
\end{array} \right)=0,
\end{eqnarray*}
we see immediately from the assumption (H2) that
\begin{eqnarray*}
a+b=b+c=c+a=0,
\end{eqnarray*}
and therefore
$$
a=b=c=0.
$$
Hence the matrix
\begin{eqnarray*}
KD=\left( \begin{array}{cccccc}
\xi_1k_{11}&\xi_2k_{12}&\xi_3k_{13} \\
\xi_1k_{21}&\xi_2k_{22}&\xi_3k_{23}\\
\xi_1k_{31}&\xi_2k_{32}&\xi_3k_{33}
\end{array} \right)
\end{eqnarray*}
is symmetric and the proof is completed.
\end{proof}

Finally, we conclude this section with the following Lemma.
\begin{lem}\label{lemma4}
Under the assumptions (H1)--(H4), there is a skew symmetric matrix
$H$ and a positive constant $c$ such that
\begin{eqnarray*}
H\Lambda - \Lambda H\ge cI - P^T\left( \begin{array}{cccccc}
0 &0  \\
0 &I_{r-1}
\end{array} \right)P.
\end{eqnarray*}
\end{lem}

\begin{proof}
let $\xi=(\xi_1, \xi_2, \cdots, \xi_r)^T$ be in the kernel of the
matrix $K$ and $A_0$ be the positive definite diagonal matrix in
Lemma \ref{lemma3}. Consider the system of equations
\begin{eqnarray*}
A_0W_t + A_0AW_x - (A_0K + K^TA_0)W = 0.
\end{eqnarray*}
From the proof of Lemma \ref{lemma3}, we know that the kernel of the
symmetric matrix $A_0K+K^TA_0$ is equal to that of $K$. By Lemma
$\ref{lemma1}$, we may assume that $\xi_i>0$ for each $i$. Thanks to
the assumption (H4), one cannot find any number $\lambda$ such that
$$
\lambda_i\xi_i=\lambda\xi_i
$$
for all $i$. Namely, $\xi=(\xi_1, \xi_2, \cdots, \xi_r)^T$ is not an
eigenvector of the coefficient matrix $\Lambda$. According to
Shizuta and Kawashima (see Theorem $1.1$ in \cite{SK} and also
Theorem $2.3$ in \cite{Y3}), there exists a positive constant $c$
and a skew symmetric matrix $H$ such that
\begin{eqnarray*}H\Lambda - \Lambda H\ge cI + A_0K + K^TA_0.
\end{eqnarray*}
Hence the lemma is proved by combining this with Lemma \ref{lemma3}.
\end{proof}

\section{Existence and Stability of steady solutions}
\setcounter{equation}{0}

In this section, we discuss the existence and stability of steady
solutions of $(\ref{pde3})$. The equations for steady solutions
$B=B(x)$ are
\begin{eqnarray}\label{bl1}
\Lambda B_x = KB.
\end{eqnarray}
Since the matrix $\Lambda$ is invertible, $W=\Lambda B$ satisfies
\begin{eqnarray}\label{bl2}
W_x = (K\Lambda ^{-1})W.
\end{eqnarray}
From the assumption (H5) and the definition of $\Lambda$, we see
that the matrix $K\Lambda ^{-1}$ satisfies the assumptions
(H1)--(H3). Thus, we can use Lemma $\ref{lemma2}$ and the Gershgorin
circle theorem to show that $K\Lambda ^{-1}$ has $(r-1)$ stable
eigenvalues and a zero-eigenvalue.

Given any boundary data $U_0(0)$, the solution for $(\ref{bl2})$ is
\begin{eqnarray*}
W(x) =\exp(K\Lambda^{-1}x)\Lambda U_0(0),
\end{eqnarray*}
and therefore,
\begin{eqnarray*}
B(x)=\Lambda ^{-1}\exp(K\Lambda^{-1}x)\Lambda U_0(0).
\end{eqnarray*}
Since $K\Lambda ^{-1}$ has $(r-1)$ stable eigenvalues and a
zero-eigenvalue, the matrix $\exp(K\Lambda^{-1}x)$ is bounded with
respect to $x\ge0$. In conclusion, given any boundary data $B(0)$,
the equations (\ref{bl1}) for steady solutions  have a unique
bounded solution $B=B(x)$.

Now we turn to discuss the stability. Set $\Phi(x,t)=U(x,t)-B(x)$.
From $(\ref{pde3})$ and $(\ref{bl1})$, we see that $\Phi$ satisfies
\begin{eqnarray}\label{eq2}
\Phi_t+\Lambda\Phi_x=K\Phi.
\end{eqnarray}
The initial and boundary data are
\begin{eqnarray}\label{value3}
\Phi(x, 0) = & U_0(x) -B(x), \nonumber\\
\Phi(0, t)= & 0.
\end{eqnarray}
From (\ref{value1}) and (\ref{bl1}), it follows that
\begin{eqnarray*}
\Lambda (U_{0x}(0)-B_x(0)) = K (U_0(0)-B(0)).
\end{eqnarray*}
Since $U_0(0)=B(0)$ and $\Lambda$ is invertible, the consistency of the initial and boundary data becomes
\begin{eqnarray}\label{value2}
\Phi_x(0,0) = \Phi(0,0) = 0.
\end{eqnarray}
Thus, our task is reduced to analyzing time-asymptotic behaviors of
the solution $\Phi$ to the IBVP $(\ref{eq2})$ together with
$(\ref{value3})$.

We start with the following local existence result. This result can
be showed by slightly modifying the proof given in Section 5 of
\cite{LY} and we omit it here.
\begin{lem}
Suppose $U_0(x)-B(x)\in H^2$ and (\ref{value2}) holds. Then there exist a positive constant
$T_*$ such that (\ref{eq2}) together with (\ref{value3}) has a unique solution
$\Phi(x,t)\in C(0,T_*;H^2)$.
Moreover, the solution satisfies the following estimate
\begin{eqnarray*}
\sup_{0\le t\le T_*}||\Phi(\cdot,t)||_{H^2}\le 2||U_0(x)-B(x)||_{H^2}.
\end{eqnarray*}
Here $T_*$ depends only on the range of $B(x)$ and any upper bound
of $||U_0(x)-B(x)||_{H^2}$.
\end{lem}
\begin{thm}
Under the assumptions (H1)--(H5), if $U_0(x)-B(x)\in H^2$ and
(\ref{value1}) holds, then (\ref{pde3}) together with (\ref{value}) has a unique global solution
$U\in C(0,\infty; H^2)$ satisfying
\begin{eqnarray*}
\lim_{t\rightarrow \infty}\sup_{x\in \mathbb{R}_{+}}|U(x,t)-B(x)|=0.
\end{eqnarray*}
\end{thm}

\begin{proof}
Let $A_0$ be the matrix given in Lemma \ref{lemma3}. We multiply
$(\ref{eq2})$ with $\Phi^TA_0$ to get
\begin{eqnarray*}
\Phi^TA_0\Phi_t+\Phi^TA_0\Lambda\Phi_x=\Phi^TA_0K\Phi.
\end{eqnarray*}
Namely,
\begin{eqnarray*}
\Phi^TA_0\Phi_t+\Phi^TA_0\Lambda\Phi_x=\frac{1}{2}\Phi^T(A_0K+K^TA_0)\Phi\le-\frac{c}{2}\Phi^TP^T\left( \begin{array}{cccccc}
0 &0  \\
0 &I_r
\end{array} \right)P\Phi,
\end{eqnarray*}
where we have used Lemma $\ref{lemma3}$. Setting
$$
V=P\Phi=(V_1, V_2, \cdots, V_r)^T
$$
and integrating the above inequality with respect
to $(x,t)\in[0, \infty)\times[0, t]$, we get
\begin{eqnarray*}
\int_0^{+\infty}{\frac{1}{2}}\Phi^TA_0\Phi(x,t) dx -
\int_0^{+\infty}{\frac{1}{2}}\Phi^T A_0\Phi(x,0) dx\\
-\int_0^t\frac{1}{2}\Phi^T A_0\Lambda \Phi(0,\tau)d\tau
+\frac{c}{2}\int_0^t\int_0^{+\infty}\sum_{i=2}^r{V_i}^2(x,\tau)dxd\tau\le0.
\end{eqnarray*}
Since $A_0$ is positive definite and $\Phi(0,t)=0$, there is a
generic constant $C$ such that
\begin{eqnarray}\label{ineq0}
\|\Phi(\cdot,t)\|_{L^2}^2 +
\sum_{i=2}^r\int_0^t{\|V_i(\cdot,\tau)\|_{L^2}}^2d\tau\le
C\|\Phi(\cdot,0)\|_{L^2}^2.
\end{eqnarray}

Next we estimate the derivatives. Differentiating the equation
$(\ref{eq2})$ with respect to $x$ gives
\begin{eqnarray}\label{De}
\Phi_{xt}+\Lambda\Phi_{xx}=K\Phi_x,\nonumber \\
\Phi_{xxt}+\Lambda\Phi_{xxx}=K\Phi_{xx}
\end{eqnarray}
From the equations (\ref{eq2}) and (\ref{De}) with the data in $(\ref{value3})$,
the boundary data for the derivative are
\begin{eqnarray*}
\Phi_x(0,t)=\Lambda^{-1}(K\Phi(0,t)-\Phi_t(0,t))=0,\\
\Phi_{xx}(0,t)=\Lambda^{-1}(K\Phi_x(0,t)-\Phi_{xt}(0,t))=0.
\end{eqnarray*}
Thus, we use the same technique shown above to estimate
$||V_x(\cdot,t)||$ and $||V_{xx}(\cdot,t)||$,
\begin{eqnarray*}
\|\Phi_x(\cdot,t)\|_{L^2}^2 +
\sum_{i=2}^r\int_0^t{\|{V_i}_x(\cdot,\tau)\|_{L^2}}^2d\tau\le
C\|\Phi_x(\cdot,0)\|_{L^2}^2,\\
\|\Phi_{xx}(\cdot,t)\|_{L^2}^2 +
\sum_{i=2}^r\int_0^t{\|{V_i}_{xx}(\cdot,\tau)\|_{L^2}}^2d\tau\le
C\|\Phi_{xx}(\cdot,0)\|_{L^2}^2.
\end{eqnarray*}
Summing these and the inequality in (\ref{ineq0}), we get
\begin{eqnarray}\label{ineq1}
\|\Phi(\cdot,t)\|_{H^2}^2 +
\sum_{i=2}^r\int_0^t\|V_i(\cdot,\tau)\|_{H^2}^2d\tau\le
C\|\Phi(\cdot,0)\|_{H^2}^2.
\end{eqnarray}

On the other hand, we multiply the equation $(\ref{eq2})$ with
$\Phi_x^TH$ to get
\begin{eqnarray}\label{eq3}
\Phi_x^TH\Phi_t + \Phi_x^TH\Lambda \Phi_x=\Phi_x^THK\Phi.
\end{eqnarray}
Here $H$ is the skew symmetric matrix in Lemma $\ref{lemma4}$. Since
\begin{eqnarray*}
\Phi_x^TH\Phi_t =
\frac{1}{2}(\Phi_x^TH\Phi)_t-\frac{1}{2}(\Phi_t^TH\Phi)_x.
\end{eqnarray*}
we integrate (\ref{eq3}) with respect to $(x,t)$ and use Lemma
\ref{lemma4} to obtain
\begin{eqnarray*}
&&c\int_0^t\int_0^{+\infty}\Phi_x^2(x,\tau)dxd\tau\le \int_0^t
\int_0^{+\infty}\sum_{i=2}^rV_{ix}^2(x,\tau)dxd\tau \\
&& + 2\int_0^t\int_0^{+\infty}\Phi_x^THK\Phi(x,\tau)dxd\tau
-\int_0^{+\infty}\Phi_x^TH\Phi(x,t)dx +
\int_0^{+\infty}\Phi_x^TH\Phi(x,0)dx.
\end{eqnarray*}
Recall that $\Phi=P^TV$ thanks to the orthogonality of $P$ and the
first column of $P^T$ is in the kernel of $K$. Thus, $K\Phi=KP^TV$
is independent of the first component of $V$. Therefore the last
inequality becomes
\begin{eqnarray*}
&&c\int_0^t\int_0^{+\infty}\Phi_x^2(x,\tau)dxd\tau\le \int_0^t
\int_0^{+\infty}\sum_{i=2}^rV_{ix}^2(x,\tau)dxd\tau\\
&& + 2\int_0^t\int_0^{+\infty}\Phi_x^THKP^TV(x,\tau)dxd\tau\\
&&-\int_0^{+\infty}\Phi_x^TH\Phi(x,t)dx +
\int_0^{+\infty}\Phi_x^TH\Phi(x,0)dx\\
&\le& \frac{c}{2}\int_0^t\|\Phi_x(\cdot,\tau)\|^2d\tau +
C\int_0^t\sum_{i=2}^r\|V_i(\cdot,\tau)\|_{H^2}^2d\tau\\
&&+C\|\Phi(\cdot,t)\|_{H^2}^2+C\|\Phi(\cdot,0)\|_{H^2}^2 .
\end{eqnarray*}
Namely,
\begin{eqnarray*}
\int_0^t\int_0^{+\infty}\Phi_x^2(x,\tau)dxd\tau \le
C\int_0^t\sum_{i=2}^r\|V_i(\cdot,\tau)\|_{H^2}^2d\tau +
C\|\Phi(\cdot,t)\|_{H^2}^2+C\|\Phi(\cdot,0)\|_{H^2}^2 .
\end{eqnarray*}
Combining this with $(\ref{ineq1})$ gives
\begin{eqnarray}\label{ineq3}
\|\Phi(\cdot,t)\|_{H^2}^2
+\int_0^t\sum_{i=2}^r\|V_i(\cdot,\tau)\|_{H^2}^2d\tau +
\int_0^t\|\Phi_x(\cdot,\tau)\|_{L^2}^2d\tau\le
C\|\Phi(\cdot,0)\|_{H^2}^2
\end{eqnarray}
for all $t\ge0$.

Furthermore, for any $t_1, t_2>0$, we have
\begin{eqnarray*}
|\|\Phi_x(\cdot,t_2)\|_{L^2}-\|\Phi_x(\cdot,t_1)\|_{L^2}|
\le\|\Phi_x(\cdot,t_2)-\Phi_x(\cdot,t_1)\|_{L^2}\\
=\|\int_{t_1}^{t_2}\Lambda\Phi_{xx}(\cdot,\tau)d\tau-
\int_{t_1}^{t_2}K\Phi_x(\cdot,\tau)d\tau\|_{L^2}\\
\le C|t_2-t_1|\max_{\tau}\|\Phi(\cdot,\tau)\|_{H^2}\le C|t_2-t_1|,
\end{eqnarray*}
where we have used (\ref{De}) in the second step. Thus, we can deduce that $\|\Phi_x(\cdot,t)\|_{L^2}\rightarrow 0 $ as
$t\rightarrow+\infty$. Moreover, from $(\ref{ineq3})$, we see that
$\|\Phi\|_{L^2}$ is bounded. Following from the celebrated
Gagliardo-Nirenberg inequality
\begin{eqnarray*}
|\Phi|_{\infty}\le\sqrt{2}\|\Phi\|_{L^2}^{\frac{1}{2}}\|\Phi_x\|_{L^2}^{\frac{1}{2}},
\end{eqnarray*}
we derive $|\Phi(\cdot,t)|_{\infty}\rightarrow 0$ as
$t\rightarrow+\infty$. Recall that $\Phi(x,t)=U(x,t)-B(x)$, and we complete the proof.
\end{proof}




\begin{thebibliography}{50}

\bibitem{Can}
P. Cancalon, \emph{Influence of temperature on slow flow in
populations of regenerating axons with different elongation
velocities}, Developmental Brain Research, 9 (1983), pp. 279-289.

\bibitem{Ca}
D. D. Carr, \emph{Global existence of solutions to
reaction-hyperbolic systems in one space dimension}, SIAM. J. Math.
Anal. {\bf 26(2)} (1995), 399--414.


\bibitem{CBF}
G. Craciun \& A. Brown \& A. Friedman, \emph{A dynamical system
model of neurofilament transport in axons}, J. Theoretical Biology
{\bf 237} (2005), 316--322.

\bibitem{DY}
A. Dressel \& W.-A. Yong, \emph{Existence of Traveling-Wave
Solutions for Hyperbolic Systems of Balance Laws}, Arch. Rational
Mech. Anal {\bf 182}, 49--75.

\bibitem{FC}
A. Friedman \& G. Craciun, \emph{Approximate travelling waves in
linear reaction-hyperbolic equations}, SIAM. J. Math. Anal. {\bf
38(3)}, (2006), 741--758.

\bibitem{FH}
A. Friedman \& B. Hu, \emph{Uniform convergence for approximate
travelling waves in linear reaction-hyperbolic systems}, Indiana
Univ. Math. J. {\bf 56 (5)} (2007), 2133--2158.

\bibitem{GB}
G.W. Gross \& L.M.Beidler, \emph{A quantitative analysis of
isotope concentration profiles and rapid transport velocities
in the C-fibers of the garfish olfactory nerve}, J.Neurobiol.,
6 (1975), pp. 213-232.

\bibitem{LY}
Hailiang Liu \& W.-A. Yong, {\em Time-Asymptotic Stability of
Boundary-Layers for a Hyperbolic Relaxation System}, Comm.PDE.
{\bf 26} (7\&8) (2001), 1323-1343.

\bibitem{Oc}
S.Ochs,  \emph{Rate offast axoplasmic transport in mammalian
nervefibers}, J. Physiol.,227 (1972), pp. 627-645.

\bibitem{RB}
M. C. Reed \& J. J. Blum, \emph{Mathematical Questions in Axonal
Transport}, In: Lectures on Mathematics in the Life Sciences, Vol.
24, 1994.

\bibitem{RVB}
M. C. Reed \& S. Venakides \&  J. J. Blum,  \emph{Approximate
travelling waves in linear reaction-hyperbolic equations}, SIAM. J.
Appl. Math. {\bf 50(1)} (1990), 167--180.

\bibitem{SK}
Yasushi Shizuta \& Shuichi Kawashima, {\em Systems of equations
of hyperbolic-parabolic type with applications to the discrete
Boltzmann equation}, Hokkaido Mathematical Journal. Vol.14, 1985, 249-275.

\bibitem{YGCGM}
G.S.Yablonsky, A.N.Gorban, D.Constales, V.Galvita, G.B.Marin
\emph{Reciprocal Relations Between Kinetic Curves
}, EPL,93 (2011).

\bibitem{Y1}
W.-A. Yong, \emph{Singular perturbations of first-order hyperbolic
systems}, Ph.D. Thesis, Universit\"{a}t Heidelberg, 1992.

\bibitem{Y2}
W.-A. Yong, {\em Basic aspects of hyperbolic relaxation systems}, in
Advances in the Theory of Shock Waves, H.~Freist\"uhler and
A.~Szepessy, eds., Progress in Nonlinear Differential Equations and
Their Applications, Vol. 47, Birkh\"auser, Boston, 2001, 259--305.

\bibitem{Y3}
W.-A. Yong, {\em Entropy and Global Existence for Hyperbolic
Banlance Laws}, Arch.Rational Mech.Anal. {\bf 172} (2004), 247-266.

\bibitem{YY}
Hao Yan \& W.-A. Yong, {\em Weak Entropy Solutions of
Nonlinear Reaction-Hyperbolic Systems for Axonal Transport},
Mathematical Models and Methods in Appl. Sci. (accepted).

\bibitem{YZ}
W.-A. Yong \& K.Zumbrun \emph{Existence of relaxation shock
profiles for hyperbolic conservation laws},
Siam.J.Appl.Math. Vol.60, 2000, no.5, 1565-1575.


\end{thebibliography}
\end{document}